\input amstex
\documentstyle{amsppt}
\input epsf.tex

\def\C{{\Bbb C}}
\def\R{{\Bbb R}}
\def\Z{{\Bbb Z}}

\def\Rp#1{\R\roman P^{#1}}
\def\Cp#1{\C\roman P^{#1}}
\def\barRP#1{\overline{\R\roman P}^{#1}}
\def\barCP#1{\overline{\C\roman P}^{#1}}
\def\conj{\mathop{\roman{conj}}\nolimits}
\def\Fix{\mathop{\roman{Fix}}\nolimits}
\def\Cl{\mathop{\roman{Cl}}\nolimits}

\def\oo{\varnothing}
\def\s{\sigma}
\def\til#1{\widetilde{#1}}
\def\pt{{\roman{pt}}}
\def\D{\Delta}

\let\tm\proclaim
\let\endtm\endproclaim
\let\rk=\remark
\let\endrk=\endremark
\let\ge\geqslant

\let\+\sqcup
\let\dsum\+
\let\-\smallsetminus
\let\d\partial

\baselineskip 20pt \magnification\magstep1 \nologo \NoRunningHeads
\NoBlackBoxes

\topmatter
\title Exotic embeddings of $\#6\Bbb R\roman{P}^2$ in the 4-sphere
\endtitle
\author S.~Finashin
\endauthor
\address Middle East Technical University,
Department of Mathematics\endgraf Ankara 06531 Turkey
\endaddress
\email  serge\@metu.edu.tr \endemail
%\date \enddate
%\thanks
%\endthanks
%\dedicatory \enddedicatory
%\keywords \endkeywords
\abstract
 We construct an infinite sequence of smooth embeddings of
$\#6\Bbb R\roman{P}^2$ in $S^4$, which are all ambient
homeomorphic, but pairwise ambient non-diffeomorphic. The double
covers of $S^4$ ramified along these surfaces form a family of the
exotic $\Bbb C\roman{P}^2\#5\overline{\Bbb C\roman{P}^2}$
constructed recently by Park, Stipsicz and Szab\'o.
 \endabstract
\endtopmatter
%\endcomment

%\rightline{\vbox{\hsize 70mm \noindent\eightit\baselineskip10pt \vskip3mm
%
% \noindent\eightrm \vskip10mm

\document
\heading\S1. Introduction\endheading

 \tm{Theorem A} For any $k\ge6$ there exists
an infinite family of smoothly embedded surfaces $F_i\subset S^4$,
$i=1,2,\dots$, homeomorphic to $F=\#k\Rp2$ {\rm(}connected sum of
$k$ copies of $\Rp2${\rm)} and with the normal Euler number
$F^2=2k-4$, such that \roster \item the pairs $(S^4,F_i)$ are all
homeomorphic; the ambient homeomorphisms can be assumed to be
diffeomorphisms in some neighborhoods of $F_i$; \item $(S^4,F_i)$
are all pairwise non-diffeomorphic.
\endroster
\endtm

Theorem A improves the result of \cite{FKV1}-\cite{FKV2}, where a
similar family of $F_i\subset S^4$ was constructed for $k=9$. Our
construction of $F_i$ is based on similar ideas and makes use of
the examples of exotic $\Cp2\#5\barCP2$ in \cite{PSS}. More
precisely, our goal can be stated as follows.

\tm{Theorem B} There exists an infinite  family of smoothly
embedded surfaces $F_i\subset S^4$, $i=1,2,\dots$ which are all
homeomorphic to a connected sum $\#6\Rp2$, such that
$\pi_1(S^4\smallsetminus F_i)=\Z/2$, and the double covers $X_i\to
S^4$ branched along $F_i$ provide an infinite family of exotic
$\Cp2\#5\barCP2$ constructed in \cite{PSS}.
\endtm

Recall the well-known diffeomorphism $\Cp2/\conj=S^4$. The image
of $\Rp2$ in the quotient-space $S^4$ represent an isotopy class
of {\it standard embeddings} of $\Rp2$ with the normal Euler
number $-2$. Another isotopy class of standard embeddings (with
the normal Euler class $2$) is obtained by reversing the
orientation of $S^4$. It is represented by the image of $\Rp2$ in
$S^4=\barCP2/\conj$, which will be denoted $\barRP2\subset S^4$. A
non-orientable surface $F\subset S^4$ will be called {\it
standard} if it splits into an ambient connected sum of such
standard embeddings, that is $F=p\Rp2\#q\barRP2$.

\subheading{Theorem B implies Theorem A} It is proven in
\cite{PSS} that $X_i$ are pairwise non-diffeomorphic. This implies
that the pairs $(S^4,F_i)$ are non-diffeomorphic to each other. It
is known moreover that for any $k\ge1$ $X_i\#k\barCP2$,
$i=1,2\dots$, remain pairwise non-diffeomorphic, which implies
that $F_i\#k\barRP2\subset S^4$ are also all ambient
non-diffeomorphic.

 The values $F_i^2=8$ of the normal Euler numbers
for $F_i$ in Theorem B can be obtained from the the signature
formula $\s(X_i)=-4=2\s(S^4)-\frac12F_i^2$. For the connected sum
$F_i\#k\barRP2$ the Euler number becomes $8+2k$.

It follows from \cite{FKV2} that the obstruction for ambient
homeomorphism of $F_i$ belongs to a finite group. This implies
that we can choose infinitely many ambient homeomorphic surfaces
$F_i$ in the infinite set of pairwise non-diffeomorphic ones,
which are  covered by exotic $\Cp2\#5\barCP2$.  \qed

\rk{Remark} In \cite{K}, it is shown that the finite ambiguity
observed in \cite{FKV2} for the exotic $10\#\Rp2$ actually
vanishes. This means that all the examples of embedded  $10\#\Rp2$
that were constructed in \cite{FKV2} are actually homeomorphic to
a standard $\Rp2\#9\barRP2$. It seems probable that the arguments
in \cite{K} can be adopted after an appropriate modification in
our case of $\#6\Rp2$. This would imply that all the examples of
$F_i$ in Theorem B are actually ambient homeomorphic to a standard
$\Rp2\#5\barRP2$.
\endrk

\subheading{Scheme of the proof of theorem B} There are several
alternative constructions of an exotic $\Cp2\#5\barCP2$ in
\cite{PSS}, and the one suitable for us is obtained by some
surgery from a rational elliptic surface, $X$, with a fiber of
type $\Bbb I_8$. The first step is a double node neighborhood knot
surgery on $X$, which yields a 4-manifold $X_K$ containing a nodal
{\it pseudo-section}. Next, $X_K$ is blown up at several points so
that we obtain a suitable chain of spheres, $C=C_1\cup\dots C_k$,
which can be rationally blowdown on the last step. Our aim is to
perform these constructions equivariantly.

In \S2, we construct a special example of a rational elliptic
surface, $X$, with a fiber $\Bbb I_8$, which is defined over
reals, and thus has an involution of the complex conjugation, $c$.
It is essential for the further constructions that the components
of the $\Bbb I_8$-fiber as well as the four remaining singular
fishtail fibers are all {\it real} (that is invariant under the
complex conjugation). In \S4, we observe that a certain
non-singular real fiber, $T$, which is constructed in \S2, can be
used for an equivariant double node knot surgery. We check that
the nodal pseudo-section $S_K$ obtained after such a surgery can
be chosen invariant under the involution. Following the
construction in \cite{PSS}, we blowup several points, which turn
out to be all real in our example of $X$. Finally we obtain as in
\cite{PSS} a chain of spheres $C=C_1\cup\dots\cup C_n$, whose
components are all $c$-invariant.

In \S3, we discuss an equivariant blowdown of such chains $C$. It
is crucial for us that the quotient $X/c$ turns out to be $S^4$
and that the quotient by the involution remains the same as we
modify the 4-manifold and the involution. So, all the involved
equivariant transformations of $X$ (knot surgery, blowup at a real
point and rational blowdown) just modify the fixed point set $F$
in the quotient $S^4$. Another crucial fact is that the
fundamental group $\pi_1(S^4\-F)=\Z/2$ is preserved under these
modifications of $F$.

More precisely, $\pi_1$ is preserved after a rational blowdown of
$C$ if we put a certain condition on $C$. This condition is
satisfied for two of the configurations proposed in Proposition
2.5 of \cite{PSS}: for $C_{79,44}$ and for $C_{89,9}$, which are
the chains
\newline
$(-2,-5,-11,-2,-2,-2,-2,-2,-2,-3,-2,-2,-3)$ and
\newline
$(-10,-11,-2,-2,-2,-2,-2,-2,-3,-2,-2,-2,-2,-2,-2,-2,-2)$.

Following the construction \cite{PSS} in the equivariant setting,
we obtain a certain 4-manifold $\widehat X$ homeomorphic but not
diffeomorphic to $\Cp2\#5\barCP2$, with an involution, $\widehat
c\:\widehat X\to\widehat X$. In the quotient $\widehat X/\widehat
c=S^4$ there is a surface $F\subset S^4$ which is the fixed point
set of $\widehat c$. Observing that $F$ is connected,
non-orientable (because $F^2=8\ne0$), and estimating its Euler
characteristic $\chi(F)=2\chi(S^4)-\chi(\widehat X)$, we deduce
that it is $\#6\Rp2$.

A sequence of the twist-knots $K_i$ that can be used for the knot
surgery on the first step of the construction (see figure 6a)
yields a sequence $\widehat X_i$ of exotic $\Cp2\#5\barCP2$ with
involutions $\widehat c_i$, and a sequence of surfaces $F_i\subset
S^4$ required for theorem B. \qed

\subheading{Acknowledgements} I thank A.Degtyarev for a discussion
about the canonical Weierstrass models of real elliptic surfaces,
which helped the author to come to the construction in \S2.

%%%%%%%%%%%%%%%%%%%%%%%%%%%%%%%%%%%%%%%%%%%%%%%%%%%%%%%%

\heading \S2. Real rational elliptic surfaces with special
singular fibers
\endheading

\subheading{2.1. Double planes ramified along quartics} Recall
that the double covering over $\Cp2$ branched along a non-singular
quartic, $A\subset\Cp2$, yields a del Pezzo surface
$X_A=\Cp2\#7\barCP2$. A pencil of lines, $L_s\subset\Cp2$,
centered at $x\in\Cp2$, is covered by an elliptic pencil,
$T_s\subset \Cp2\#7\barCP2$, whose singular fibers correspond to
the lines tangent to $A$. Blowing up the pull-back of $x$ in
$\Cp2\#7\barCP2$, we obtain an elliptic fibration, $p\:X\to\Cp1$.

Assume that the quartic $A$ has a singular point, $y\in A$, of the
type $\Bbb A_n$ (by definition at such a singularity $p$ has a
local model $(z_1,z_2)\mapsto z_1^{n+1}+z_2^2$), and the basepoint
$x$ is generic. Then we obtain an elliptic fibration with a
singular fiber $T_0$ of the type $\Bbb I_{n+1}$ (in Kodaira's
classification). The fiber $T_0$
 corresponds to the line $L_0\subset\Cp2$ passing through $x$ and
 $y$. If in addition $L_0$ is tangent to $A$ at the basepoint $x\in
A$, then the fiber $T_0$ is of the type $\Bbb I_{n+3}$.

For instance, we obtain $\Bbb I_8$-fiber $T_0$ if we choose a
quartic $A=A_1\cup A_3$ that splits into a cubic, $A_3$, and a
real line, $A_1$, tangent to $A_3$ at its inflection point, $y$
(see Figure 1a). The corresponding line $L_0$ should pass through
$y$ and be tangent to $A_3$ at some other point, $x$, which will
be the center of the pencil of the lines $L_s$.

\subheading{2.2. Construction of a special real elliptic
fibration} Consider the double covering $q\:X_A\to\Cp2$ ramified
along a degree $2n$ curve $A\subset\Cp2$ having an equation $f=0$
with real coefficients. The complex conjugation in $\Cp2$ can be
lifted to $X_A$ in two ways.
 Namely, there are two real algebraic models of $X_A$
defined by a weighted homogeneous equation $y^2=\pm
f(x_0,x_1,x_2)$ in a quasi-projective space $P_{1,1,1,n}$ with the
coordinates $x_0,x_1,x_2,y$ of weights $1,1,1,n$. The
corresponding involutions $c_\pm\:X_A\to X_A$, induced from the
complex conjugation in $P_{1,1,1,n}$, have fixed point
 sets $\Fix(c_\pm|_{X_A})$ which are projected by $q$
to the regions $\Rp2_{\pm f}=\{x\in\Rp2\,|\,\pm f(x)\ge0\}$
bounded by the curve $A_\R=A\cap\Rp2$.

In our example of the real quartic $A=A_1\cup A_3$, we choose the
region $\Rp2_f$ as is shown on Figure 1c) and the corresponding
involution $c_A\:X_A\to X_A$ whose fixed point set is
$q^{-1}(\Rp2_f)$.

\midinsert \epsfbox{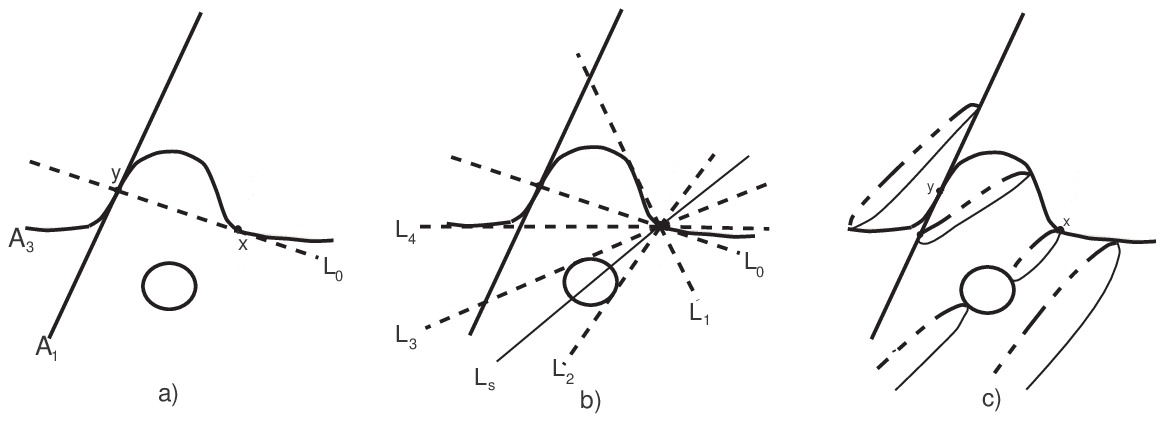} \vskip-4mm\botcaption{Figure
1}\endcaption \eightpoint \noindent  a) Quartic $A=A_1\cup A_3$
and the basepoint $x$. b) The tangent lines $L_i$, $i=1,2,3,4$ of
the pencil and line $L_s$ corresponding to the fiber $T=T_s$. c)
The real locus of the double plane $X_A$.
\endinsert

Blowing up the singularity and the two infinitely near base-points
of $X_A$, we obtain a {\it real elliptic fibration},
$p\:X\to\Cp1$, endowed with an involution $c\:X\to X$ commuting
with $p$. Let $F=\Fix(c)$ denote its fixed point set.

\tm{Lemma 1} The real elliptic fibration $p\:X\to\Cp1$ constructed
above has the following properties. \roster\item $X$ contains a
real singular fiber $T_0=C_1\cup\dots\cup C_8$ of the type $\Bbb
I_8$, whose components $C_i$ are $c$-invariant. \item $X$ contains
4 other  $c$-invariant singular fibers, $T_i$, $i=1,2,3,4$, which
are ordinary fishtails. \item The elliptic fibration $p\:X\to\Cp1$
admits a  $c$-invariant section.
 \item $X$ can be blowdown to $\Cp2$, so that each of the
nine consecutively contracted $(-1)$-curves is real.
 This transforms the non-singular real fibers $F_s$ to non-singular real cubics.
 \item $X/c=S^4$.
 \item
If the singular fibers $T_i$, $i=0,1,2,3,4$ are the pre-images of
the tangent lines $L_i$ on Figure 1b, then a non-singular real
fiber $T=T_s$ chosen between $T_2$ and $T_3$ has real locus,
$T\cap F$, formed by two connected components, like is shown on
Figure 2c. The two vanishing curves in $T$, which are contracted
as $T$ is degenerated into the singular fibers $T_2$ and $T_3$,
are isotopic. A vanishing curve from this isotopy class can be
chosen $c$-invariant, and so that $c$ reverses its orientation.
 \item The complement $F\-(F\cap T)$ is
connected.
\endroster
\endtm

\demo{Proof} The real form of the singularity ${\Bbb A}_5$
involved in our construction has local model $x^6-y^2+z^2=0$, and
$c$-invariance of the components $C_i$ is verified by its
analysis. The fibers $T_i$ in are $c$-invariant because the
tangent lines $L_i$ are real.

To justify (3), we will present 6 real sections. The first one is
the exceptional curve, $E_x\subset X$, which appears after the
second of the two infinitely near blowups at the basepoint $x\in
X_A$, as we construct $X$. Another section is the proper
transform, $A_1^*\subset X$, of the line $A_1$. The proper
transform of a real line $L'\subset\Cp2$ passing through $y$ and
tangent to $A_3$ (see Figure 2a) splits in two components, $L'_1$,
$L_2'$, which are both sections of $p$. Another tangent line $L''$
shown on Figure 2a give similarly components $L''_1$ and $L_2''$.

\midinsert \epsfbox{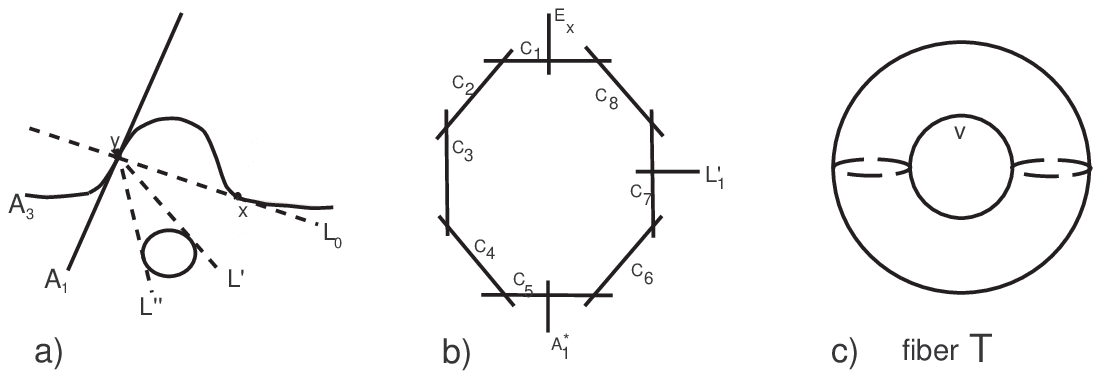} \vskip-4mm\botcaption{Figure 2}
\endcaption
\noindent \eightpoint a) Tangent lines $L'$ and $L''$. b) Fiber
$T_0$ and 3 disjoint real sections $E_x$, $A_1^*$, and $L_1'$. c)
The fixed point set of $c$ dividing the fiber $T$, and the
vanishing curve $v$.\newline
\endinsert

Let us choose the cyclic order of the components $C_i$ of $T_0$,
and suppose that $C_1$ intersects $E_x$.
 Then $A_1^*$ intersects $C_5$, and
the curves $L_1'$ and $L_1''$ intersect $C_7$ and $C_3$. We may
suppose that $L_1'$ intersect $C_7$, like it is shown on Figure 2.
If we blowdown consecutively $E_x$, $C_1$, $C_2$, $C_3$, $C_4$,
$A_1^*$, $L_1'$, $C_7$, and $C_6$, then we obtain $\Cp2$, which
proves (4). (Note that the remaining components, $C_5$ and $C_8$,
will represent a line and a conic in $\Cp2$ obtained after $9$
blowdowns.)

We can deduce (5) from (4) using that $\Cp2/\conj=S^4$, which
implies that a blowup at a real point effects to the quotient as a
connected sum with $S^4=\barCP2/\conj$, and thus, topologically
does not change the quotient.

Inspecting Figure 1c we observe that the fixed point set $T\cap F$
of the complex conjugation acting on $T$ has two connected
components, as it is shown on Figure 2c. This can be understood
from Figure 1c. The type of the vanishing curves on $T$ is
determined by the types of the real critical points of the
projection $F\to \Rp1$ (restriction of $p$). The critical points
in the fibers $T_2$ and $T_3$ have both index one, which implies
that the homology class of the corresponding vanishing curves
belong to the $(-1)$-eigenspace of $c_*$ in $H_1(T)$. This
determines these vanishing curves up to isotopy, thus, proving
(6).

Property (7) is clear from Figure 1c), if we take into account
connectedness of the real locus of the corresponding singularity
$x^6-y^2+z^2=0$ after its resolution.\qed\enddemo

%%%%%%%%%%%%%%%%%%%%%%%%%%%%%%%%%%%%%%%%%%%%%%%%%%%%%

\heading \S3. Equivariant Rational Blowdown\endheading

\subheading{3.1. Rational blow-down surface surgery} Let $X$
denote a smooth 4-manifold with {\it a chain of spheres}
$C=C_1\cup\dots\cup C_k\subset X$, which intersect each other
consecutively and transversely, so that their dual weighted graph
is a chain-tree sketched on Figure 3. A regular neighborhood,
$N(C)$, of $C$ is a plumbing 4-manifold, $P_C$, corresponding to
this weighted graph.

Certain chains $C$ can be {\it rationally blowdown}, that is we
can remove $N(C)=P_C$ from $X$ and replace it by some rational
homology ball, $Q_C$, with the same boundary $\d Q_C=\d N(C)$.
This gives a new 4-manifold $\widehat X=X'\cup Q_C$, where
$X'=X\-N(C)$.

\midinsert
\epsfbox{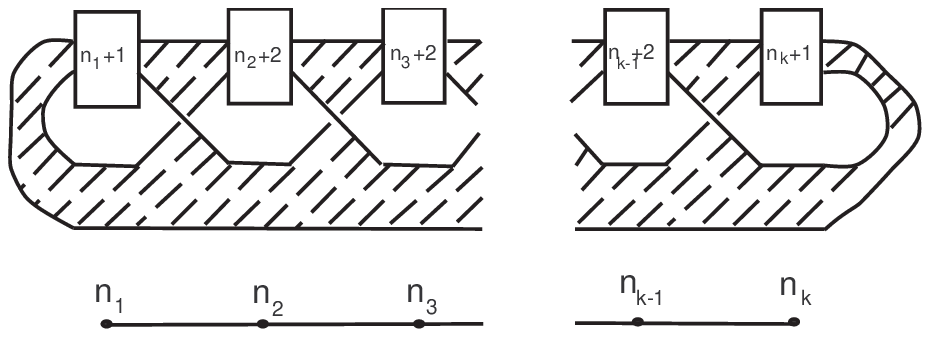} \vskip-4mm\botcaption{Figure
3}\endcaption \eightpoint \noindent  The plumbing surface $F_C$
described by a chain-tree can be presented as the span-surface of
a two-bridge knot diagram. The numbers in the boxes count the
half-twists.
\endinsert

It is well known and easy to see that $P_C$ can be described as
the double cover over $D^4$ branched along a surface, $F_C$,
obtained by plumbing of the unknotted bands, $F_{n_i}\subset D^4$,
 $i=1,\dots,k$, where $n_i$ stands for the framing of the band (number of
its half-twists which is taken with sign ``$-$'' in the case of
left-hand half-twists). Such a plumbing surface can be sketched as
is shown on Figure 3.

As it is observed in \cite{FS1}, $Q_C$ is the double cover of
$D^4$ branched along another surface, $R_C\subset D^4$, bounded by
the same link as $F_C$, $L_C=\d R_C=\d F_C$.
 More details about surface $R_C$ can be extracted from \cite{CH},
and we will only mention that $R_C$ is obtained by connecting a
pair of disjoint unknotted discs, $D_1\+ D_2$, via a ribbon in
$S^3$ and then pushing in inside $D^4$ the interior of the
surface. This ribbon connects either $\d D_1$ with $\d D_2$, or
the boundary of one of the discs to itself, say of $\d D^1$. In
the first case we obtain a knotted disc in $D^4$. In the second
case we obtain a disjoint union of a disc $D_2$ with a M\"obius
band, as it is shown on Figure 4b in the simplest example. (One
can show that the band cannot be orientable, because $\d F_C$
cannot have 3 boundary components.)

\midinsert \epsfbox{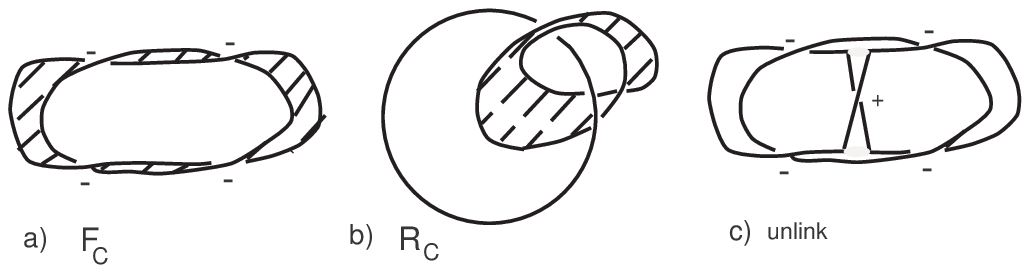} \vskip-4mm\botcaption{Figure
4}\endcaption \eightpoint \noindent a) A band $F_C=F_{-4}$ with
the four negative half-twists (marked with the signs ``$-$''). b)
$R_C$ (a disc and a band) bounded by the same link $L_C$. c) $L_C$
differs from an unlink by a ribbon move.
\endinsert

\tm{Lemma 2} Consider a rational blowdown of a chain of spheres
$C\subset X$, which yields $\widehat X=X'\cup R_C$. Assume that
$X$ is endowed with an orientation preserving involution $c\:X\to
X$, which keeps each of the components, $C_i\subset C$ invariant,
and reverses its orientation, so that $\Fix(c)\cap C_i\ne\oo$.
Then, if $N(C)$ is chosen $c$-invariant, the rational blowdown can
be made equivariantly, which yields an involution $\widehat
c\:\widehat X\to\widehat X$.

Such a blowdown gives the same quotient $Y=X/c=\widehat X/\widehat
c$. The fixed point sets $F=\Fix c$ and $\widehat F=\Fix(\widehat
c)$ descended to $Y$ give the same locus $F'=F\cap Y'=\widehat
F\cap Y'$, inside $Y'=X'/c$. Furthermore, $N(C)/c=D^4$, and $F\cap
D^4$ is isotopic to the plumbing surface $F_C$.
 The piece of surface $\widehat F\cap D^4$ is isotopic to
the surface $R_C$.
\endtm

\demo{Proof} Under these assumptions, $c|_{N(C)}$ is equivalent to
the deck transformation of the double branched covering $N(C)\to
D^4$. The involution $\widehat c$ just extends the involution
 $c|_{X\-N(C)}$ to $Q_C$ as the deck transformation of the
branched covering $Q_C\to D^4$. \qed\enddemo

We say that $(\widehat X,\widehat c)$ is obtained by an {\it
equivariant rational blowdown} from $(X,c)$.

\subheading{3.2. The characteristic sub-configuration} It is not
difficult to see that the number of components of the link $L_C=\d
F_C$ is determined by the parity of numbers $C_i^2$, namely, $L_C$
has one component if the intersection matrix $(C_i\circ C_j)$ is
non-singular modulo $2$, and has two component if singular.

We say that the union of some of the components $C_i$ forms a {\it
characteristic sub-configuration}, $W\subset C$,  if the
fundamental class $[W]\in H_2(C;\Z/2)$ is a Wu element of the
intersection form $(C_i\circ C_j)$, that is $C_i^2=C_i\circ
W\mod2$ for all $i\in\{1,\dots,k\}$.
 It is easy to see that
the characteristic sub-configuration $W$ is unique if the matrix
$(C_i\circ C_j)$ is non-singular modulo $2$ (has odd determinant).
If this matrix is $\mod2$ singular, then there are two
characteristic sub-configurations, $W$ and $W'$, whose sum gives
the non-trivial element of the null-space of $(C_i\circ C_j)$ in
$H_2(C;\Z/2)$ (it is easy to check that this null-space has
dimension at most $1$).

\rk{Remark} It is simple to determine $W$ using an orientation of
the link diagram of $L_C$. Namely, $W$ contains $C_i$ if and only
if the opposite sides of the band $F_{n_i}$ are co-directed,
 like is shown on Figure 5. This follows from that $W\cap F$
realizes the first Stiefel-Whitney class $w_1(F_C)$.\endrk

\subheading{3.3. Commutativity of $\pi_1$ after rational
blowdowns}
 Suppose that a sphere $C_0\subset X$ extends the chain $C$ to a longer chain
$\til C=C_0\cup C_1\cup\dots\cup C_n$. This means that $C_0$ is a
$c$-invariant sphere (like the other $C_i$) which intersects $C_1$
transversely at a single point and does not intersect $C_i$, if
$i>1$.

\tm{Lemma 3} Assume that \roster\item the link $L(C)$ is a knot,
 \item a characteristic sub-configuration
$W\subset C$ is not characteristic for $\til C$.
\endroster
Then $\pi_1(Y\- \widehat F)=\pi_1(Y\- F)$.
\endtm

\rk{Remark} The second assumption of the lemma means that $C_1$ is
not included into $W$ if $n_0=C_0^2$ is odd, and $C_1\subset W$ if
$n_0$ is even.
\endrk

\demo{Proof} Applying the Van Kampen theorem, we see that
 $\pi_1(Y\- F)=G'*_{G_L} G_C$, where $G'=\pi_1(Y'\- F')$, $G_C=\pi_1(D^4\-
 F_C)$,
and $G_L=\pi_1(S^3\- L_C)$. Similarly, $\pi_1(Y\- \widehat
F)=G'*_{G_L}\widehat G$, where $\widehat G=\pi_1(D^4\- R_C)$.

The plan of the proof is to observe that the inclusion
homomorphisms $G_L\to G_C$, $G_L\to \widehat G$ are epimorphisms,
and that their kernels, $K$ and $\widehat K$, vanish under the
inclusion homomorphism $G_L\to G'$. This implies that the
homomorphisms $G'\to G'*_{G_L} G_C$ and $G'\to G'*_{G_L}\widehat
G$ are isomorphisms.

First of all, note that the upper Wirtinger presentation for the
link $L_C$ diagram shown on Figure 3 implies that its group $G_L$
is generated by two elements $a,b$, presented by the loops around
the overpasses $\ell_a$ and $\ell_b$ shown of  Figure 5.

The homomorphism $G_L\to G_C$ is epimorphic because $G_C$ is a
cyclic group (since $F_C$ is a connected span-surface for $L_C$,
whose interior is pushed out from $S^3$ inside $D^4$). Inspecting
the homology, we see that $G_C=\Z$ is obtained from $G_L$ by
adding the relation $a=b$, in the case of orientable surface
$F_C$. If $F_C$ is non-orientable, then $G_C=\Z/2$ is obtained
from $G_L$ by adding two relations: $a=b$ and $a=b^{-1}$. These
relations generate $K$.

The homomorphism $G_L\to \widehat G$ is also epimorphic, because
$R_C$ is a ribbon-surface. The kernel $\widehat K$ is contained in
the commutator subgroup $[G_L,G_L]$, which is the kernel of the
product map $G_L\to\widehat G\to H_1(D^4\-R_C)=H_1(S^3\-L_C)=\Z$,
where the latter two equalities are due to our assumption that
$L_C$ is connected, and so $R_C$ is a disc. Thus, $[G_L,G_L]$ is
generated by the relation $a=b$ if overpasses $\ell_a$ and
$\ell_b$ on Figure 5 inherit co-directed orientation from $L_C$,
and by the relation $a=b^{-1}$ if these overpasses inherit
opposite orientations.

\midinsert \epsfbox{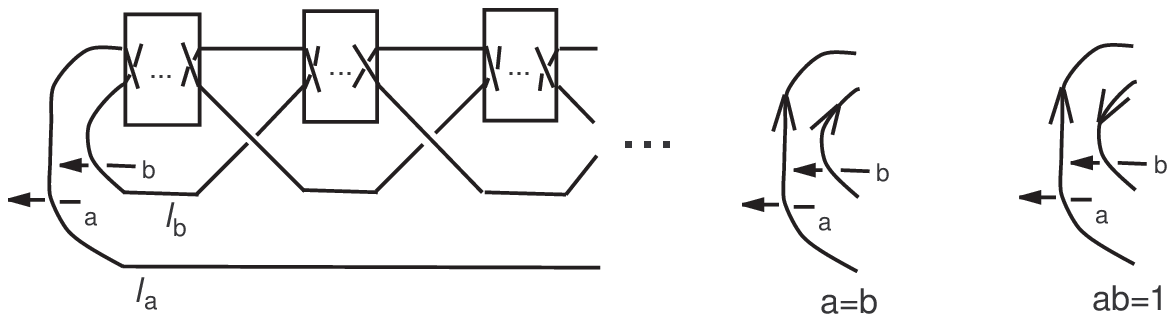} \vskip-4mm\botcaption{Figure
5}
\endcaption{\eightpoint \noindent Overpasses $\ell_a$,
$\ell_b$ and the corresponding generators $a$ and $b$ of
$\pi_1(S^3\-L_C)$. The case of co-directed and oppositely directed
overpasses $\ell_a$, $\ell_b$, with the corresponding relations
between $a$ and $b$ (after adding the commutativity relation
$ab=ba$) }
\endinsert

Showing that the images of $a,b$ under the inclusion homomorphisms
$G_L\to G'$ (for which we keep the same notation $a,b$) satisfy
the both relations $a=b$ and $a=b^{-1}$, will complete the proof.

One of these relations comes from a regular neighborhood $N(D_0)$
of the disc $D_0=C_0/c$ in $Y$. Note that $H=N(D_0)\cap Y'$ is a
4-ball containing an unknotted disc $F_H=F\cap H$, so that
$\pi_1(H\-F_H)=\Z$. The common piece of boundary of $H$ and $D^4$
is a 3-ball, which intersects  $F$ along a pair of arcs,
$\ell_a\cup\ell_b$. It is not difficult to see that in
$\pi_1(H\-F_H)$ we obtain the relation $a=b$ if $n_0$ is even, and
$a=b^{-1}$ if odd. With this relation, group $G_L$ becomes
abelian, and we obtain another relation (which depends on the
orientation of $\ell_a$, $\ell_b$ induced from $L_C$, as was
explained). Under the second assumption of our lemma, these two
relations are different, that is the both relations $a=b$ and
$a=b^{-1}$ are satisfied in $\pi_1(Y\-\widehat F)$.
 \qed\enddemo

%%%%%%%%%%%%%%%%%%%%%%%%%%%%%%%%%%%%%%%%%%%%%%%%%%%%%%%%%%%%%%%%%%%%%%%%

\heading \S4. Equivariant version of the Fintushel-Stern double
node knot surgery
\endheading

\subheading{4.1. Equivariant knot surgery} The $4$-dimensional
knot surgery consists in removing from a $4$-manifold $X$ a
trivialized tubular neighborhood $N(T)=T\times D^2$ of a torus
$T\subset X$, and replacing it by $S^1\times C(K)$, where $C(K)$
is a knot complement (see \cite{FS1}). It is supposed that the
gluing map $S^1\times\d C(K)\to T\times \d D^2$ identifies
 a longitude $\pt\times\ell\subset S^1\times C(K)$
with a meridian of $T$, $m_T=\pt\times\d D^2\subset\d N(T)$, which
yields a 4-manifold $X_K$, homologically equivalent to $X$.

In the equivariant version of this construction, we suppose that
$X$ is endowed with an orientation preserving involution, $c$,
which keeps invariant the torus $T$ as well as its neighborhood
$N(T)$.
We say that a trivialization $N(T)=T\times D^2$ of a tubular
neighborhood, $N(T)$, of $T$ is {\it equivariant} if the action of
$c$ on $N(T)$ is presented as the direct product of $c|_T$ and the
complex conjugation in $D^2\subset\C$. Note that equivariant
trivializability is equivalent to existence of a projection
$N(T)\to D^2$ which commutes with $c|_{N(T)}$ and the complex
conjugation in $D^2$. In the case of our interest, $T$ is a real
non-singular fiber in a real elliptic fibration and thus admits
such an equivariantly trivializable neighborhood.

Let us assume in addition that $c|_T$ reverses orientation of $T$
and has two-component fixed point set, $F\cap T$ (see Figure 2c).
In this case the quotient $\Cal A=T/c$ is an annulus and in the
coordinates defined by some diffeomorphism $T=S^1\times S^1$ the
action of $c$ looks as $(z_0,z_1)\mapsto(z_0,\bar z_1)$. Thus for
an appropriate diffeomorphism $N(T)=S^1\times S^1\times D^2$, this
action looks as $(z_0,z_1,z_2)\mapsto(z_0,\bar z_1,\bar z_2)$.

From a knot $K\subset S^3$ we require that it has an axis of
symmetry, and intersects this axis at a pair of points. It will be
convenient to choose the complex conjugation, $\conj$, in
$S^3\subset\C^2$ as such a symmetry, so that the axis is
$S^1_\R=S^3\cap\R^2$.
Such a knot $K$ admits an equivariant tubular neighborhood, $N(K)$
in which $c$ acts as $(z_1,z_2)\mapsto(\bar z_1,\bar z_2)$, with
respect to a trivialization $S^1\times D^2=N(K)$. We can choose
such a trivialization to be {\it null-framed}, which means that a
longitude $\ell_K=S^1\times\pt$ is null-homologous in the knot
complement $C(K)=\Cl(S^3\-N(K))$.

We can glue $S^1\times C(K)$ to $X\-N(T)$ via an equivariant
gluing map $g\:S^1\times\d C(K)\to \d N(T)$. Using the coordinates
$(z_0,z_1,z_2)$ in $\d C(K)=\d N(K)$ and $\d N(T)$ from the above
trivializations of $N(K)$ and $N(T)$, we define map $g$ as
$(z_0,z_1,z_2)\mapsto (z_0,z_2,z_1)$. Such an equivariant knot
surgery yields a 4-manifold $X_K$ endowed with an involution,
$c_K\:X_K\to X_K$.

\subheading{4.2. The tangle surgery in the quotient-spaces} The
quotient $N(K)/\conj$ is a 3-ball, which can be viewed as a
regular neighborhood, $N(\frak s_K)$, of the arc $\frak
s_K=K/\conj$ in $S^3=S^3/\conj$. Thus, $B_K=C(K)/\conj$ is also a
3-ball, complementary to $N(\frak s_K)$.
 The unknot $S^1_\R\subset S^3/\conj$ splits
into a trivial tangle $\frak t=S^1_\R\cap N(s_K)$ and a
non-trivial one, $\frak t_K=S^1_\R\cap B_K$ (see Figure 6c).

\rk{Example 1} The twist-knot $K=K_n$, which will be used in our
construction, admits a $\conj$-invariant presentation, as it is
sketched on Figure 6a. Figures 6d--6e present the corresponding
tangle splitting $S^1_\R=\frak t\cup \frak t_K\subset S^3$.
\endrk

\vskip-2mm \midinsert \epsfbox{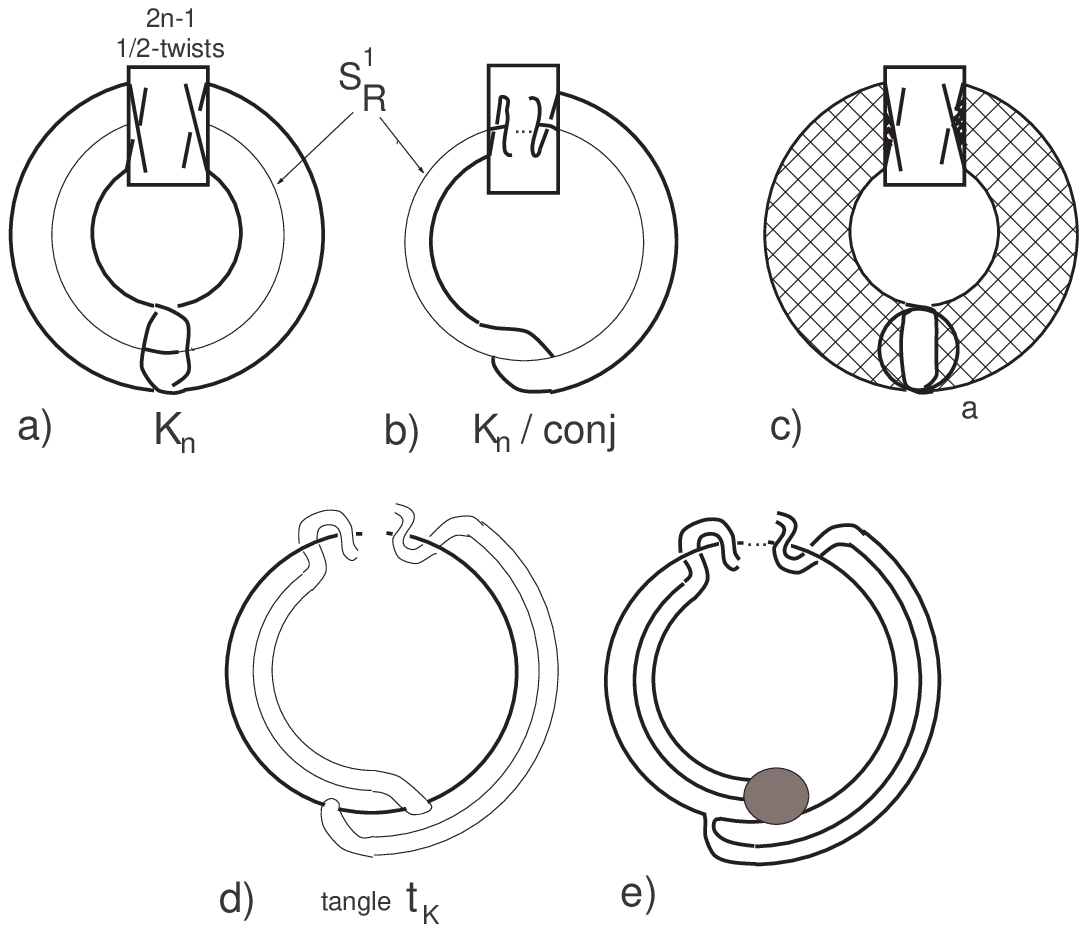}\vskip-4mm
\botcaption{Figure 6}
\endcaption\eightpoint \noindent  a) The twist-knot $K=K_n$ with the axis of symmetry
$S^1_\R$. b) The arc $\frak s_K=K/\conj$. c) Seifert surface
$S_K^\circ$ of genus 1 bounded by $K$. It contains
$\conj$-invariant
 curve $a$, which bounds a $\conj$-invariant disc $D_a$.
 d) The ball
$N(\frak s_K)$ and tangle $\frak t_K$ in its complement. e)$\frak
t_K$ after deformation of $N(\frak s_K)$ (the ball shaded on the
figure).
\endinsert

The quotient-space $X_K/c_K$ is obtained from $X/c$ by removing a
regular neighborhood, $N=N(\Cal A)$, of the annulus $\Cal A=T/c$,
which can be viewed as $N=S^1\times N(\frak s_K)=S^1\times B^3$,
and replacing it by $S^1\times B_K=S^1\times B^3$. Such a surgery
does not change the topological type of a 4-manifold, so we can
identify the both quotients, $Y=X/c=X_K/c_K$.

 The branching locus $F_K$ of the double covering $X_K\to Y$
is obtained from $F$ after replacing $F_N=F\cap N=S^1\times \frak
t$ by $S^1\times \frak t_K$ inside $N=S^1\times B^3$. Note that
the components of $\frak t_K$ and of $\frak t$ connect the same
pairs of their common endpoints. We denote these four endpoints
$p_1^\pm,p_2^\pm$, and assume that $p_i^+$ is connected with
$p_i^-$.
 Moreover, the both tangles must have
{\it the same framing}. This means by definition that the kernel
of the inclusion homomorphism $H_1(S^2\-\d \frak t)\to H_1(D^3\-
\frak t)$ is the same as for $H_1(S^2\-\d \frak t_K)\to H_1(D^3\-
\frak t_K)$.

Such kind of surgery will be called {\it tangle surgery of
$F\subset Y$ along an annulus membrane $\Cal A$. It can be applied
to any surface $F$ in a 4-manifold $Y$ under the assumption that
the annulus membrane $\Cal A$ with the boundary on this surface is
{\it null-framed}. This
 means by definition that for some trivialization
$N=S^1\times D^3$ of its regular neighborhood, $N=N(\Cal A)$, the
part of surface $F\cap N$ is identified with $S^1\times \frak t$,
and $\Cal A$ is identified with $S^1\times \frak s$, where $\frak
s$ is a line segment connecting the midpoints of the components of
$\frak t$ (see Figure 7a). The following lemma summarizes our
observations.

\tm{Lemma 4} An equivariant knot surgery on $(X,c)$ along a
$c$-invariant torus $T\subset X$ gives $(X_K,c_K)$ with the same
quotient-space $Y=X_K/c_K=X/c$. The fixed point set $F_K\subset Y$
of $c_K$ is obtained from the fixed point set $F$ by the tangle
surgery along the annulus membrane $\Cal A=T/c$. \qed\endtm

\subheading{4.3. Commutativity of $\pi_1$ throughout the knot
surgery}
 \tm{Lemma 5}
Assume that $F\subset Y$ is a surface in a 4-manifold and $\Cal A$
is a null-framed annulus membrane on $F$ such that $F\-\d\Cal A$
is connected. Assume that $F_K$ is obtained from $F$ by applying
the tangle surgery along $\Cal A$ with respect to the $\frak t_K$,
where $K=K_n$ is the twist-knot from Example 1. Assume furthermore
that the group $\pi_1(Y\-(F\cup\Cal A))$ is abelian. Then
$\pi_1(Y\- F_K)$ is also abelian and isomorphic to $\pi_1(Y\- F)$.
\endtm

\midinsert \epsfbox{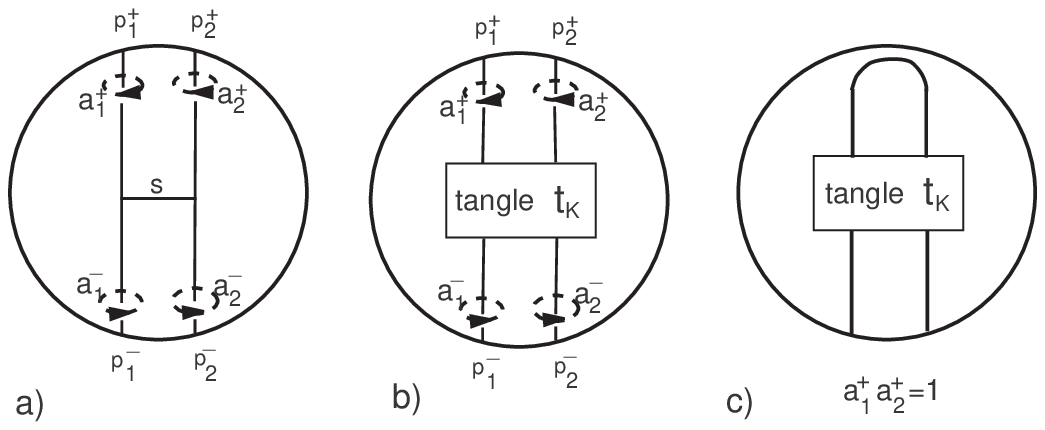}\vskip-4mm
\botcaption{Figure 7}
\endcaption\eightpoint \noindent a) Trivial tangle $\frak t$
with the connecting line segment $\frak s$. The generators
$a_i^\pm$ of $\pi_1(S^2\-\d \frak t)$. b) The result of a tangle
surgery. c) Adding the relation $a_1^+a_2^+=1$ to the group of
tangle $\frak t_K$ effects like connecting together the points
$p_1^+$ and $p_2^+$.
\endinsert

\demo{Proof} Let $Y'=\Cl(Y\-N)$ and $F'=F\cap Y'$,
$F_{N,K}=F_K\cap N$.
 By the Van Kampen theorem,
$\pi_1(Y\-F_K)=G'*_HG_N$, where $G'=\pi_1(Y'\-F')$, $H=\pi_1(\d
N\-\d F_N)$, and $G_N=\pi_1(N\-F_{N,K})$. Note that $N\-(\Cal
A\cup F_N)$ can be deformation retracted to its boundary $\d N\-\d
F_N$, which implies that group $G'$ is abelian, by the assumption
on $\pi_1(Y\-(F\cup\Cal A)$. Note that $G'*_HG_N=G'*_H(G_N/K)$,
where $K$ is the image in $G_N$ of the kernel of the homomorphism
$H\to G'$. We will show that $G_N/K$ is an abelian group and the
product homomorphism $H\to G_N\to G_N/K$ is epimorphic. This
implies that $G'*_HG_N$ is a quotient of $G'$ and thus is also
abelian.

Note that $H=\Z\times\pi_1(S^2\-\d\frak t)$, where the second
factor is a free group or rank 3. It is convenient to present this
free group by 4 generators, $a_1^\pm$, $a_2^\pm$, satisfying the
relation $a_1^+a_1^-a_2^-a_2^+=1$. These generators correspond to
the loops around the tangle endpoints, $p_1^\pm,p_2^\pm\in S^2$,
in the positive direction, see Figure 7a.

Let us fix some element $a\in G'$ presented by a loop around
$F\-\d\Cal A$. Commutativity of $G'$ and connectedness of
$F\-\d\Cal A$ imply that the inclusion homomorphism $H\to G'$
sends each of $a_i^\pm$ either to $a$, or to $a^{-1}$ (depending
on the topology of the boundary $\d\Cal A$ as an oriented curve in
$F$). Such a relation, $a_i^\pm=a$ or $a_i^\pm=a^{-1}$, is
inherited by the quotient-group $G_N/K$. To complete proof of the
lemma it is enough to show that by adding this relation we
transform group $\pi_1(D^3\-\frak t_K)$ into a cyclic group with a
generator $a$ (since the factor $\Z$ in
$G_N=\Z\times\pi_1(D^3\-\frak t_K)$ lies in the center and comes
from the corresponding factor in
$H=\Z\times\pi_1(S^2\-\{p_1^+,p_1^-,p_2^+,p_2^-\})$).

We will present two arguments. The first one can be applied to any
knot $K$ admissible for an equivariant knot surgery, but it works
only if we have a relation $a_1^+a_2^+=1$, or $a_1^-a_2^-=1$. Note
that if we connect the endpoints $p_1^+$ and $p_2^+$ as is shown
on Figure 7c, we modify the group $\pi_1(B^3\-\frak t_K)$ by
adding a relation $a_1^+a_2^+=1$. In the case of tangles $\frak
t_K$ constructed from $\conj$-invariant knots $K$, this
modification transforms the tangle into an unknotted arc in $D^3$.
Thus, the group $\pi_1(D^3\-\frak t_K)$ becomes cyclic and
generated by any of the elements $a_i^\mp$. The case of adding
relation $a_1^-a_2^-=1$ is analogous.

Our second argument is specific for the twist-knot $K=K_n$, but
can be applied in the case of relation $a_1^+=a_2^+$ or
$a_1^-=a_2^-$ (as well as in the case of relation $a_1^\pm
a_2^\pm=1$ considered before). First, we observe that the upper
Wirtinger presentation gives 5 generators for $\pi_1(D^3\- \frak
t_K)$, namely $a_i^\pm$, $i=1,2$, and one more generator $b$ shown
on Figure 8.

The two strands of the tangle $\frak t_K$ with the origins at the
points $p_1^+$ and $p_2^+$ pass together several times under
$S^1_\R$. These underpasses separate the consecutive overpasses on
the first strand which yield elements $a_1^+, b^{-1}a_1^+b,\dots$,
which are all conjugate to $a_1^+$. The similar overpasses on the
second strand give elements $a_2^+,b^{-1}a_2^+b,\dots$, which are
conjugate to $a_2^+$ by the same sequence of elements. In the end
of the sequence, we obtain elements $b^{-1}=x^{-1}a_1^+x$ and
$(a_1^-)^{-1}=x^{-1}a_2^+x$, which are conjugate to $a_1$ and
$a_2$ via the same element $x\in\pi_1(D^3\- \frak t_K)$.
 So, a relation $a_1^+=a_2^+$ implies that $b=a_1^-$, whereas
$a_1^+a_2^+=1$ implies $ba_1^-=1$. In any case, generator $b$ can
be eliminated, and after adding two more relations to $\pi_1(D^3\-
\frak t_K)$, namely  $a_1^-=a_1^+$ (or $a_1^-=(a_1^+)^{-1}$) and
$a_2^-=a_1^+$ (or $a_2^-=(a_1^+)^{-1}$, we obtain a cyclic group.

Finally, we can observe that surface $F_K$ is homologically
equivalent to $F$, and thus $H_1(Y\-F_K)=H_1(Y\-F)$. Since the
group $\pi_1(Y\-F)$ is obviously abelian due to the assumption of
the lemma, we obtain an isomorphism $\pi_1(Y\-F_K)=\pi_1(Y\-F)$.
\qed\enddemo

\vskip-2mm\midinsert
\centerline{\epsfbox{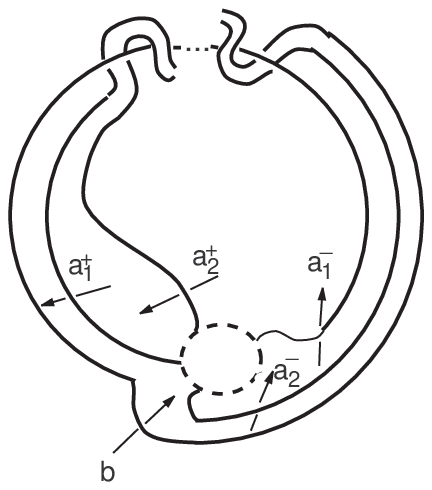}}\vskip-5mm
\botcaption{Figure 8}
\endcaption \eightpoint \noindent
The tangle group becomes abelian after adding the relations
$a_1^+=a_2^+$, $a_1^-=(a_1^+)^{\pm1}$, and $a_2^-=(a_1^+)^{\pm1}$.
\endinsert

\tm{Lemma 6} Let $(X,c)$ be the real elliptic surface constructed
in \S2, and $T$ the real fiber from Lemma 1(6). Then membrane
$\Cal A=T/c$ satisfies the assumptions of Lemma 5, and thus
$\pi_1(S^4\-F_K)=\Z/2$ (here $S^4=X/c$ by Lemma 1(6)).\endtm

\demo{Proof} Connectedness of $F\-\Cal A$ is observed in Lemma
1(4).

The group $\pi_1(S^4\-(F\cup\Cal A))$ was shown to be cyclic in
\cite{FKV2}, \S4, under the assumption that $T\subset X$ is
obtained from a real non-singular cubic in $\Cp2$ by blowing up
the base-points of a real pencil of cubics. This is so in our
case, as follows from property (6) of Lemma 1.
 \qed\enddemo

\subheading{4.4. The equivariant double node surgery}
 To justify that a pseudo-section
$S_K\subset X_K$ can be chosen $c$-invariant we recall first its
construction in \cite{FS2}. Consider a disc $\D_1\subset \Cp1$
which contains inside precisely two critical values
$s_+,s_-\subset\D_1$ of an elliptic Lefschetz fibration
$p\:X\to\Cp1$. Assume moreover that the corresponding two
vanishing curves in a non-singular fiber, $T_s$, $s\in \D_1$, are
isotopic. Let $\D\subset\D_1$ denote a smaller disc not containing
points $s_\pm$, and $U=p^{-1}(\D)$. Consider a section $S\subset
X$ of $p$. Its restriction over $\D$ is the disc $\til{\D}=S\cap
U$. Since the gluing map in the definition of the knot surgery
which yields $X_K$ may be changed by an isotopy, we can make the
boundary $\d \til{\D}$ match with the boundary of a Seifert
surface $S_K^\circ\subset C(K)\times\pt\subset X_K$ and obtain a
closed surface $S_K^*=(S\-\til{\D})\cup S_K^\circ$ in $X_K$.
 If $K=K_n$ (the twist-knot on Figure 6a),
then $S_K^*$ is a torus which has a certain disc membrane
$D^*_a\subset X_K$ bounded by curve $a=\d D^*_a=D_a^*\cap S_K^*$
and having self-intersection $(D_a^*)^2=-1$ (relative to the
boundary on the surface $S_K^*$).
 The torus $S_K^*$ can be deformed into a fishtail $S_K\subset
 X_K$, as we pinch curve $a\subset S_K^*$ along disc $D_a^*$.
 The local topology of $S_K$ near its singular point is like near
an algebraic double point, and the embedded surface $S_K\subset
X_K$ is topologically equivalent to a rational curve with a single
node and self-intersection $(S_K)^2=-1$.
% (the latter is because $S_K$ is homologous to $S$).

To construct disc $D^*_a$, we first take a disc $D_a\subset S^3$
bounded by $a$, so that $D_a$ intersects $K$ at a pair of points
(see Figure 6c). Disc $D_a$ punctured at these points is embedded
in $pt\times C_K\subset S^1\times C_K\subset X_K$. The boundary of
the punctures are curves in two different fibers of $p$, and by
definition of our knot surgery these curves are vanishing
(corresponding to the singular values $s_\pm$) and so can be
filled by the discs centered at the nodes of the singular fibers
over $s_\pm$.

\tm{Lemma 7} Consider the real elliptic surface $(X,c)$
constructed in \S2. Let $K=K_n$ be the twist-knot embedded
$\conj$-invariantly in $S^3$, as is shown on Figure 6a. Assume
that $(X_K,c_K)$ is obtained from $(X,c)$ by an equivariant knot
surgery along the non-singular fiber $T=T_s$ specified in Lemma
1(6). Then, the pseudo-section $S_K$ can be chosen
$c_K$-invariant.
\endtm

\demo{Proof} By Lemma 1(3), we can suppose that section $S$ is
$c$-invariant. We consider a disc $\D\ni s$ which is invariant
under the complex conjugation in $\Cp1$, and denote by $r_\pm$ the
endpoints of the interval $\D\cap\Rp1$.
 In our example of real elliptic surface in \S2, we have a pair of
real critical values, $s_\pm$, whose fibers $T_2=p^{-1}(s_-)$,
$T_3=p^{-1}(s_+)$ can be used for a double node knot surgery, as
it follows from Lemma 1(6).
 The curve $S\cap\d U$ is a
$\conj$-invariant longitude of the knot $K$ in the boundary of
$C_K\cong\pt\times C_K$. This longitude spans a $\conj$-invariant
Seifert surface $S_K^\circ\subset C_K$, as it is shown on Figure
6c. This implies that torus $S_K^*$ can be chosen $c_K$-invariant.
Furthermore, we can choose the disc $D_a^*$ to be $c_K$-invariant
as well. Namely, we choose first a $\conj$-invariant disc $D_a$
(see Figure 6c) and make $\conj$-invariant punctures around the
two intersection points $D_a\cap K$, which are both real. The
boundary of such punctures are $c$-invariant curves in the two
fibers $p^{-1}(r_\pm)\subset\d U$, namely, the vanishing curves of
the critical values $s_\pm$.
 If the discs filling these curves are chosen
$c$-invariant, then the disc $D_a^*$ becomes $c$-invariant as
well.

Finally, note that there is a $c_K$-equivariant deformation of
$S_K$, which contracts disc $D_a^*$ and degenerates the torus
$S_K$ into $S_K^*$. Its construction goes like in the
non-equivariant case: we deform $s_K$ using a flow of a vector
field tangent to $D_a^*$. To obtain an equivariant deformation,
this field should be chosen $c$-invariant.
 \qed\enddemo

%%%%%%%%%%%%%%%%%%%%%%%%%%%%%%%%%%%%%%%%%%%%%%%%%%%%

\enddocument